\documentclass[10pt,reqno]{amsproc}
\newtheorem{theorem}{\bf{Theorem}}[section] 
\newtheorem{lemma}[theorem]{\bf{Lemma}}     

\newtheorem{corollary}[theorem]{\bf{Corollary}}

\newtheorem{definition}[theorem]{\bf{Definition}}
\newtheorem{remark}[theorem]{\bf{Remark}}

\newcommand*{\blue}[1]{{\color{blue}#1}}
\newcommand*{\red}[1]{{\color{red}#1}}

\usepackage{color}

\title[OPERATOR VALUED MAPS ON HILBERT $C^*$-MODULES]
 {OPERATOR VALUED MAPS ON HILBERT $C^*$-MODULES} 
\author[M. B. Asadi]{Mohammad B. Asadi}
\address{School of Mathematics, Statistics and Computer Science,
 College of Science, University of Tehran, Tehran,
  Iran, and \\
   School of Mathematics, Institute for Research in Fundamental Sciences (IPM),
P.O. Box: 19395-5746, Tehran, Iran}

\email{mb.asadi@khayam.ut.ac.ir}

\author[R. Behmani]{Reza Behmani}
\address{\noindent Department of Mathematics, Kharazmi University, 50, Taleghani Ave.,15618, Tehran Iran}

\email{reza.behmani@gmail.com}

\author[A. R. Medghalchi]{Ali R. Medghalchi}
\address{\noindent Department of Mathematics, Kharazmi University, 50, Taleghani Ave.,15618, Tehran IRAN.}

\email{a\_medghalchi@khu.ac.ir}%

\author[H. Nikpey]{Hamed Nikpey}
\address{Department of Mathematics, Shahid Rajaei Teacher Training University, Tehran 16785-136, Iran}
\email{hamednikpey@gmail.com}

\subjclass[2010]{Primary: 46L08, Secondary: 46L07.}
\keywords{Hilbert $C^*$-modules, Completely positive maps,
$\varphi$-maps}

\begin{document}
\maketitle

\begin{abstract}
We provide a characterization for operator valued completely
bounded linear maps on Hilbert $C^*$-modules in terms of
$\varphi$-maps. Also, we show that for every operator valued
completely positive map $\varphi$ on a $C^*$-algebra
$\mathcal{A}$, there is a unique (up to multiplication by a
unitary operator) non-degenerate $\varphi$-map on each Hilbert
$\mathcal{A}$-module.

 \end{abstract}




\section{INTRODUCTION} 

The study of $\varphi$-maps on Hilbert $C^*$-modules has increased
significantly during recent decades. In this context, several
concepts such as representation theory of Hilbert $C^*$-modules,
dilation theory of $\varphi$-maps and CP-extendable maps were
studied
 (\cite{Arambašic, Asadi, ABN, BRS, Dey-Tervidi, Heo-Ji, J1, J2, skeide, Skeide-Sumesh}).
Therefore it becomes natural to concentrate on $\varphi$-maps as
important maps on Hilbert $C^*$-modules. To confirm this
statement, we show that an operator valued map on a Hilbert
$C^*$-module is completely bounded if and only if it can be
decomposed  to a bounded operator and a $\varphi$-map, for a
completely positive map $\varphi$ on the underlying $C^*$-algebra
of the Hilbert $C^*$-module.

Moreover, for a given operator valued completely positive map
$\varphi$ on a $C^*$-algebra $\mathcal{A}$ and its minimal
Stinespring dilation $\pi$, we construct a $\varphi$-map and a
$\pi$-representation for each Hilbert $\mathcal{A}$-module
$\mathcal{E}$ and also we show that every  non-degenerate
$\varphi$-map ( non-degenerate $\pi$-representation) on
$\mathcal{E}$ is a unitary operator multiple of the above
constructed $\varphi$-map ($\pi$-representation).

We denote Hilbert spaces by $\mathcal{H,K,L}$. The set of all
bounded operators between Hilbert spaces $\mathcal{H,K}$ is
denoted by $\mathcal{B(H,K)},$ and
$\mathcal{B(H)}:=\mathcal{B(H,H)}.$

  Assume  $\mathcal{E}$ is a
Hilbert $C^*$-module over a unital $C^*$-algebras $\mathcal{A}.$
The linking $C^*$-algebra of $\mathcal{E}$ is denoted by
$\mathcal{L}(\mathcal{E})$ and defined as
$\mathcal{L}(\mathcal{E}):=\{\begin{bmatrix}
T & x \\
y^* & a
\end{bmatrix} : \ a\in\mathcal{A}, T\in\mathbb{K}(\mathcal{E}), x,y\in\mathcal{E}\},$
where $\mathbb{K}(\mathcal{E})$ is the set of compact operators on
$\mathcal{E}$.
 Also, for arbitrary given maps
$\rho:\mathcal{A}\rightarrow\mathcal{B(H)}$,
$\sigma:\mathbb{K}(\mathcal{E})\rightarrow\mathcal{B(K)}$ and
 $\Psi:\mathcal{E}\rightarrow\mathcal{B}(\mathcal{H},\mathcal{K})$,
 the map $\begin{bmatrix}
T & x \\
y^* & a
\end{bmatrix} \mapsto \begin{bmatrix}
\sigma(T) & \Psi(x) \\
\Psi(y)^* & \rho(a)
\end{bmatrix}$ from $\mathcal{L}(\mathcal{E})$ into $\mathcal{B}(\mathcal{K} \oplus
\mathcal{H})$ is denoted by $\begin{bmatrix}
\sigma & \Psi \\
\Psi^* & \rho
\end{bmatrix}$.

If $\varphi:\mathcal{A}\rightarrow\mathcal{B(H)}$ is a completely
positive map and $\Phi:\mathcal{E}\rightarrow\mathcal{B(H,K)}$
 a map, then we say

(1) $\Phi$ is \textit{non-degenerate}, if
$[\Phi(\mathcal{E})\mathcal{H}]=\mathcal{K}.$

(2) $\Phi$ is a \textit{$\varphi$-map}, if $\Phi(x)^*\Phi(y)=
\varphi(\langle x,y\rangle)$, for all $x,y\in\mathcal{E}$.

(3) $\Phi$ is a \textit{representation} (or \textit{$\rho$-representation}),
 if there is a $*$-representation $\rho:\mathcal{A}\to\mathcal{B(H)}$ such that $\Phi$ is a $\rho$-map.

(4) $\Phi$ is a \textit{completely semi-$\varphi$-map}, if
$\Phi_n(x)^*\Phi_n(x)\leq\varphi_n(\langle x,x\rangle)$ for every $n
\in \mathbb{N}$ and $x\in\mathbb{M}_n(\mathcal{E}).$

(5) $\Phi$ is a \textit{CP-extendable map}, if there exist
completely positive maps
$\phi_1:\mathbb{K}(\mathcal{E})\to\mathcal{B(K)}$ and
$\phi_2:\mathcal{A}\rightarrow\mathcal{B(H)}$  such that
$\begin{bmatrix}
\phi_1 & \Phi \\
\Phi^* & \phi_2
\end{bmatrix}:\mathcal{L(E)}\to\mathcal{B(K\oplus H)},$
is a completely positive map.

(6) $\Phi$ is \textit{dilatable} if there is a   representation
$\Psi:\mathcal{E}\to\mathcal{B}(\mathcal{H'},\mathcal{K'})$ and
bounded operators $V:\mathcal{H}\to\mathcal{H'}$ and
$W:\mathcal{K}\to\mathcal{K'}$ such that $\Phi(x)=W^*\Psi(x)V,$
for every $x\in\mathcal{E}.$

Positive definite kernels are a non-linear version of completely
positive maps which are older than their linear counterpart (see
\cite{Ars, Gs, kol, Pa1}) for more details). A positive definite
kernel on a set $X$ is a two variables function $\phi:X\times X
\to\mathcal{B(H)}$, where $\mathcal{H}$ is a Hilbert space, such
that for every choice of n elements in $X$ such as
$\{x_{i}\}_{i=1}^{n}$, $[\phi(x_{i},x_{j})] \in
\mathbb{M}_{n}(\mathcal{B}(\mathcal{H}))_{+}.$ From now on we use
PD kernel to abbreviate positive definite kernel.

For a given PD kernel $\phi:X\times X \to\mathcal{B(H)}$ there is
a standard way to construct another Hilbert space $\mathcal{K}$
such that $\phi$ is decomposed into more tractable functions from
$X$ into $\mathcal{B(H,K)}$ \cite{Gs,kol}.

 \begin{definition}\label{definition0}
Let $X$ be a non-empty set and $\phi:X\times X\to\mathcal{B(H)}$
be a PD kernel. A Kolmogorov decomposition pair for $\phi$ is a
pair $(\nu,\mathcal{K})$ consists of a Hilbert space $\mathcal{K}$
a map \linebreak $\nu:X\to\mathcal{B}(\mathcal{H},\mathcal{K})$
such that $\phi(x,y)=\nu(x)^{\ast}\nu(y).$ A Kolmogorov
decomposition pair is called minimal when
$[\nu(X)\mathcal{H}]=\mathcal{K}.$
 \end{definition}
The existence of the Kolmogorov decomposition pair for a PD kernel
is a well-known result:
\begin{theorem}
 Let $\phi:X\times X \to\mathcal{B(H)}$ be a PD kernel,
  then there is a Hilbert space $\mathcal{K}$ and a map $\nu:X\to\mathcal{B(H,K)}$ such that
 $\phi(x,y)=\nu(x)^*\nu(y)$,
 for all $x,y\in X.$
\end{theorem}


 \begin{remark}\label{remark0}
Minimal Kolmogorov decomposition pairs of $\phi$ are unique up
to unitary equivalence. That is, if $(\nu,\mathcal{K})$ is a
minimal Kolmogorov decomposition of $\phi$ and
$(\upsilon,\mathcal{L})$ is an arbitrary Kolmogorov
decomposition pair of $\phi$, then there is a unique isometry
$V:\mathcal{K}\to\mathcal{L}$ such that
 $V\nu(x)=\upsilon(x).$

To every map $\Phi:X\to\mathcal{B(H,K)},$ one can associate a PD
kernel $\Lambda_\Phi:X\times X\to\mathcal{B(H)}$ by
$\Lambda_\Phi(x,y)=\Phi(x)^*\Phi(y),$ which has $\Phi$ as its
Kolmogorov  decomposition. Also, a completely positive map
$\varphi:\mathcal{A}\to\mathcal{B(H)}$ induces a PD kernel
$\tilde{\varphi}:\mathcal{E}\times\mathcal{E}\to\mathcal{B(H)}$ by
$\tilde{\varphi}(x,y)=\varphi(\langle x,y\rangle)$, on every
Hilbert $\mathcal{A}$-module $\mathcal{E}.$
 \end{remark}

\section{Main Theorems}

The following theorem says that each operator valued completely
bounded map on a Hilbert $C^*$-module is an operator multiple of
some $\varphi$-map. We mention that a similar discussion can be
found in \cite[Section 3]{Skeide-Sumesh}. In fact, the main idea
of the proof is to use the fact that the space $B(K\oplus H)$ is
injective in the category of operator systems.
\begin{theorem} \label{t1}
Let $\mathcal{E}$ be a right Hilbert $C^*$-module over a unital
$C^*$-algebra $\mathcal{A}$ and
$\Phi:\mathcal{E}\to\mathcal{B(H,K)}$ be a map. The following
statements are equivalent

(i) $\Phi$ is a completely bounded linear map.

(ii) There is a completely positive map
$\varphi:\mathcal{A}\to\mathcal{B(H)}$ and a $\varphi$-map
$\Gamma:\mathcal{E}\to\mathcal{B(H,L)}$ and a bounded operator
$S:\mathcal{L}\to\mathcal{K}$ such that $\Phi(x)=S\Gamma(x)$, for
all $x\in\mathcal{E}.$
\end{theorem}

The following lemma provides a representation theorem for
completely positive maps on $C^*$-algebras, in term of maps on
Hilbert $C^*$-modules.
\begin{lemma} \label{l1}
Let $\mathcal{E}$ be a right Hilbert $C^*$-module over a unital
$C^*$-algebra $\mathcal{A}$,
$\varphi:\mathcal{A}\to\mathcal{B}(\mathcal{H})$  a completely
positive map and $(\pi,\mathcal{K},V)$ be the minimal Stinespring
dilation triple of $\varphi$. Then,
 there exists a triple
$((\Phi_{\varphi},\mathcal{H}_{\varphi}),(\Psi_{\pi},\mathcal{K}_{\pi}),W_{\varphi})$
consists of Hilbert spaces $\mathcal{H}_{\varphi}$ and
$\mathcal{K}_{\pi},$ a unitary operator
$W_{\varphi}:\mathcal{H}_{\varphi}\to\mathcal{K}_{\pi},$ a
non-degenerate $\varphi$-map
$\Phi_{\varphi}:\mathcal{E}\to\mathcal{B}(\mathcal{H},\mathcal{H}_{\varphi})$
and a non-degenerate $\pi$-representation
$\Psi_{\pi}:\mathcal{E}\to\mathcal{B}(\mathcal{K},\mathcal{K}_{\pi})$
such that $\Phi_{\varphi}(\cdot)=W_{\varphi}^*\Psi_{\pi}(\cdot)V.$
\end{lemma}

 Now, we summarize some
results about $\varphi$-maps on Hilbert $C^*$-modules. In fact, in
the following theorem, the part (i) is the same as
Bhat-Ramesh-Sumesh's theorem \cite[Theorem 2.1]{BRS} and  also
says that for every completely positive map on a $C^*$-algebra
$\mathcal{A}$, such as $\varphi:\mathcal A\to\mathcal{B(H)}$,
there is a unique (up to multiplication by a unitary operator)
non-degenerate $\varphi$-map on each Hilbert $\mathcal{A}$-module.
 The part (ii) strengthen
\cite[Theorem 3.4]{ABN} and characterizes completely
semi-$\varphi$-maps as operator multiple of $\varphi$-maps.
 Also,
the part (iii) is a similar result to (i) and finally, the part
(iv) exposes the relation between every pair of $\varphi$-maps and
$\pi$-representations on a same Hilbert $C^*$-module.

\begin{theorem} \label{t2}
With the notations of the above lemma, one has

 (i) a map
$\Phi:\mathcal{E}\to\mathcal{B}(\mathcal{H},\mathcal{H'})$ is a
(non-degenerate) $\varphi$-map if and only if
 there exist a (unitary)
isometry $S_{\Phi}:\mathcal{H}_{\varphi}\to\mathcal{H'},$ and a
(unitary) coisometry $W:\mathcal{H'}\to\mathcal{K}_{\pi}$ such
that $S_{\Phi}\Phi_{\varphi}=\Phi$ and
$\Phi(\cdot)=W^*\Psi_{\pi}(\cdot)V;$

(ii) a map
$\Phi:\mathcal{E}\to\mathcal{B}(\mathcal{H},\mathcal{H'})$ is a
(non-degenerate) completely semi-$\varphi$-map if and only if
there exist a (dense range) contraction
$S:\mathcal{H}_{\varphi}\to\mathcal{H'},$ and a (injective)
contraction $W:\mathcal{H'}\to\mathcal{K}_{\pi}$ such that $S
\Phi_{\varphi}=\Phi$ and $\Phi(\cdot)=W^*\Psi_{\pi}(\cdot)V;$

(iii) a map
$\Psi:\mathcal{E}\to\mathcal{B}(\mathcal{K},\mathcal{K'})$ is a
(non-degenerate) $\pi$-representation if and only if
 there exists an
(unitary) isometry $S_{\Psi}:\mathcal{K}_{\pi} \to \mathcal{K'}$
such that $\Psi(\cdot)=S_{\Psi}\Psi_{\pi}(\cdot);$

(iv) if $\Psi:\mathcal{E}\to\mathcal{B}(\mathcal{K},\mathcal{K'})$
is a
 $\pi$-representation and $\Phi:\mathcal{E}\to\mathcal{B}(\mathcal{H},\mathcal{H'})$
is a $\varphi$-map, then there exists a  partial isometry
$W:\mathcal{H'}\to\mathcal{K}'$ such that
$\Phi(\cdot)=W^*\Psi(\cdot)V$.
  Moreover $W$ is unitary when $\Phi$ and $\Psi$ are non-degenerate.
\end{theorem}

Completely semi-$\varphi$-maps introduced in \cite{ABN} as
generalizations of $\varphi$-maps. By the above theorem, every
completely semi-$\varphi$-map can be dilated to a representation
of the Hilbert $C^*$-module and therefore it is a linear map.
Also, CP-extendable maps introduced in \cite{Skeide-Sumesh} and
the authors in \cite[Theorem 4.2]{ABN} showed that each operator
valued map on a Hilbert $C^*$-module is dilatable if and only if
it is CP-extendable. Therefore, we have the following result.

\begin{corollary} \label{c1}
Let $\mathcal{E}$ be a right Hilbert $C^*$-module over a unital
$C^*$-algebra $\mathcal{A}$ and
$\Phi:\mathcal{E}\to\mathcal{B(H,K)}$ be a  map. The following
statements are equivalent

(i) $\Phi$ is a completely bounded linear map.

(ii)  There is a completely positive map $\varphi: \mathcal{A} \to
\mathcal{B(H)}$ and a $\varphi$-map $\Gamma: \mathcal{E} \to
\mathcal{B(H,L)}$ and a bounded operator $S: \mathcal{L} \to
\mathcal{K}$ such that $\Phi(x) = S\Gamma(x)$, for all $x \in
\mathcal{E}$.

(iii) There is a completely positive map $\psi:\mathcal{A} \to
\mathcal{B(H)}$ such that $\Phi$ is a completely semi-$\psi$-map.

(iv) $\Phi$ is dilatable.

(v) $\Phi$ is CP-extendable.

\end{corollary}

The following corollary is a well known theorem on completely
bounded maps on $C^*$-algebras \cite{Pau}. However, we can
conclude it as a special case of the above result, since each
$C^*$-algebra is a right Hilbert $C^*$-module over itself and also
$\mathbb{K}(\mathcal{A})\cong\mathcal{A},$
$\mathbb{M}_2(\mathcal{A})\cong\mathcal{L(A)}.$

\begin{corollary}
Let $\mathcal{A}$ be a unital $C^*$-algebra. If
$\psi:\mathcal{A}\to\mathcal{B(H)}$ is a completely bounded map,
there exist completely positive maps
$\phi_i:\mathcal{A}\to\mathcal{B(H)}$, $i=1,2$, such that the map
$\begin{bmatrix}
\phi_1 & \psi \\
\psi^* & \phi_2
\end{bmatrix}:\mathbb{M}_2(\mathcal{A})\to\mathcal{B(H\oplus H)}$
is completely positive.
\end{corollary}

\section{Proofs}

\blue{\textbf{Proof of Theorem \ref{t1}:}}

\red{$(i)\Rightarrow(ii):$} Assume $\Phi$ is completely bounded.
Since $\mathcal{B(K\oplus H)}=\begin{bmatrix}
   \mathcal{B(K)} & \mathcal{B(H,K)} \\
   \mathcal{B(K,H)} & \mathcal{B(H)}
   \end{bmatrix},$
   we can consider $\Phi$ as a map from $\mathcal{E}$ into $\mathcal{B(K\oplus H)}.$

   Let $\mathcal{L}_1(\mathcal{E}):=\{\begin{bmatrix}
T & x \\
y^* & a
\end{bmatrix} : \ a\in\mathcal{A}, T\in \mathbb{K}_{1}(\mathcal{E}):=\mathbb{K}(\mathcal{E})+\mathbb{C}I_{\mathcal{E}}, x,y\in\mathcal{E}\},$
be the unitization of the linking $C^*$-algebra of $\mathcal{E},$
then $\varphi$ can be extended to a completely bounded map
$\Psi:\mathcal{L}_1(\mathcal{E})\to\mathcal{B(K\oplus H)}$ by
Wittstock's extension theorem.
    Then there is a $*$-representation
$\pi:\mathcal{L}_1(\mathcal{E})\to\mathcal{B(L)}$ and bounded
operators $V_i:K\oplus H\to\mathcal{L},$ $i=1,2$ such that
$\Psi(X)=V_1^*\pi(X)V_2$ for every
$X\in\mathcal{L}_{1}(\mathcal{E}).$ Using \cite[Proposition
3.1]{Arambašic} $\mathcal{L}$ decomposes to
$\mathcal{L}_2\oplus\mathcal{L}_1$ for two orthogonal closed
subspaces $\mathcal{L}_1$ and $\mathcal{L}_2$ and there exist
$*$-representations
$\rho:\mathcal{A}\to\mathcal{B}(\mathcal{L}_1),$
$\sigma:\mathbb{K}_1(\mathcal{E})\to\mathcal{B}(\mathcal{L}_2)$
and a  $\sigma$-$\rho$-representation
$\Gamma_0:\mathcal{E}\to\mathcal{B}(\mathcal{L}_1,\mathcal{L}_2)$
such that
  $\pi=\begin{bmatrix}
\sigma & \Gamma_0 \\
\Gamma_{0}^* & \rho
\end{bmatrix}: \mathcal{L}_{1}(\mathcal{E}) \rightarrow \mathcal{B}(\mathcal{L}_2 \oplus \mathcal{L}_1).$

 Since $\Psi$ is an extension of $\Phi,$ and the operators $V_i,$ $i=1,2$ has the matrix decompositions
  $V_i=\begin{bmatrix}
S_{i,1} & S_{i,2} \\
S_{i,3} & S_{i,4}
\end{bmatrix}\in\mathcal{B(K\oplus H},\mathcal{L}_2\oplus\mathcal{L}_1),$ $i=1,2,$ one has
$$\begin{bmatrix}
0 & \Phi(x) \\
0 & 0
\end{bmatrix}=\Psi(\begin{bmatrix}
0 & e \\
0 & 0
\end{bmatrix})=\begin{bmatrix}
S_{1,1} & S_{1,2} \\
S_{1,3} & S_{1,4}
\end{bmatrix}^*\begin{bmatrix}
\sigma(0) & \Gamma_{0}(x) \\
\Gamma_{0}(0)^* & \rho(0)
\end{bmatrix}\begin{bmatrix}
S_{2,1} & S_{2,2} \\
S_{2,3} & S_{2,4}
\end{bmatrix},$$
for every $x\in\mathcal{E}.$
 Thus
$\Phi(x)=S_{1,1}^*\Gamma_{0}(x)S_{2,4}$  for every
$x\in\mathcal{E}.$

Now, if we set $\varphi(\cdot)= S_{2,4}^*\rho(\cdot)S_{2,4}$ and
$\Gamma(\cdot)=\Gamma_{0}(\cdot)S_{2,4}$, then $\varphi$ is a
completely positive map, $\Gamma$ is a $\varphi$-map and
$\Phi(\cdot)=S_{1,1}^*\Gamma(\cdot)$.

\red{$(ii)\Rightarrow(i):$} Let $\Phi(\cdot)=S\Gamma(\cdot)$.
Let $[x_{ij}] \in M_n(\mathcal{E})$, then
\begin{align*}
\Phi_{n}([x_{ij}])^*
\Phi_{n}([x_{ij}]) &= [\Gamma(x_{ji})^*S^*][S\Gamma(x_{ij})]=
 [\Gamma(x_{ji})^*]diag(S^*, \cdot\cdot\cdot, S^*) diag(S,
\cdot\cdot\cdot, S)[\Gamma(x_{ij})] \\&
 \leq ||S||^2 [\Gamma(x_{ji})^*][\Gamma(x_{ij})]
 = ||S||^2 \varphi_n(\langle
[x_{ij}],[x_{ij}]\rangle)
\end{align*}
Therefore,
$$\|\Phi_{n}([x_{ij}])\|^2
=\| \Phi_{n}([x_{ij}])^*
\Phi_{n}([x_{ij}]) \|
\leq
 ||S||^2 \| \varphi\|_{cb} \|[x_{ij}]\|^2$$
 and then $\Phi$ is a completely bounded map.

 $\hspace{15cm} \square$

\vspace{0.5 cm}

\blue{\textbf{Proof of Lemma \ref{l1}:}}

Define
$\tilde{\varphi}:\mathcal{E}\times\mathcal{E}\to\mathcal{B}(\mathcal{H})$
by
$$\tilde{\varphi}(x,y):=\varphi(\langle x,y \rangle_{\mathcal{A}})$$
for all $x,y \in\mathcal{E}.$ Note that the $\mathcal{A}$-valued
inner-product on $\mathcal{E},$ $\langle \cdot, \cdot
\rangle_{\mathcal{A}}:\mathcal{E}\times\mathcal{E}\to\mathcal{A}$
is a PD kernel and $\varphi$ is a completely positive map on
$\mathcal{A},$ therefore   $\tilde{\varphi}$ is a PD kernel on
$\mathcal{E}.$ There is a (unique) minimal Kolmogorov
decomposition  $(\Phi_{\varphi},\mathcal{H}_{\varphi})$ for
$\tilde{\varphi},$ consists of a Hilbert space
$\mathcal{H}_{\varphi}$ and a map
$\Phi_{\varphi}:\mathcal{E}\to\mathcal{B}(\mathcal{H},\mathcal{H}_{\varphi})$
such that the linear span of
$\Phi_{\varphi}(\mathcal{E})\mathcal{H}$ is a dense subspace of
$\mathcal{H}_{\varphi}$ and
$\tilde{\varphi}(x,y)=\Phi_{\varphi}(x)^*\Phi_{\varphi}(y)$ for
all $x,y \in\mathcal{E}.$ Thus $\Phi_{\varphi}$ is a
non-degenerate $\varphi$-map from $\mathcal{E}$ into
$\mathcal{B}(\mathcal{H},\mathcal{H}_{\varphi}).$

  Similarly, define $\tilde{\pi}:\mathcal{E}\times\mathcal{E}\to\mathcal{B}(\mathcal{K})$ by
$\tilde{\pi}(x,y):=\pi(\langle x,y \rangle_{\mathcal{A}})$, for
all $x,y \in\mathcal{E}.$ A similar argument like above lines
implies the existence of a (unique) minimal Kolmogorov
decomposition pair $(\Psi_{\pi},\mathcal{K}_{\pi})$ for
$\tilde{\pi},$ consists of a Hilbert space $\mathcal{K}_{\pi}$ and
a map
$\Psi_{\pi}:\mathcal{E}\to\mathcal{B}(\mathcal{H},\mathcal{K}_{\pi})$
such that the linear span of $\Psi_{\pi}(\mathcal{E})\mathcal{H}$
is a dense subspace of $\mathcal{K}_{\pi}$ and
$\tilde{\pi}(x,y)=\Psi_{\pi}(x)^*\Psi_{\pi}(y)$ for all
$x,y\in\mathcal{E}.$ Thus $\Psi_{\pi}$ is a non-degenerate
$\pi$-map from $\mathcal{E}$ into
$\mathcal{B}(\mathcal{H},\mathcal{K}_{\pi}).$

Since $(\pi,V,\mathcal{K})$ is a dilation triple for $\varphi$ and
$\Phi_{\varphi}$ is a $\varphi$-map, for every $x,y
\in\mathcal{E},$ we have
\begin{equation}\label{eq2}\Phi_{\varphi}(x)^*\Phi_{\varphi}(y)=\varphi(\langle
x,y \rangle_{\mathcal{A}})=V^*\pi(\langle
x,y\rangle_{\mathcal{A}})V=V^*\Psi_{\pi}(x)^*\Psi_{\pi}(y)V.\end{equation}
The above equation implies that for every
$x_{1},...,x_{n}\in\mathcal{E}$ and
$h_{1},...,h_{n}\in\mathcal{H}$
$$\|\sum_{i=1}^{n}\Phi_{\varphi}(x_{i})h_{i}\|_{\mathcal{H}_{\varphi}}=\|\sum_{i=1}^{n}\Psi_{\pi}(x_{i})Vh_{i}\|_{\mathcal{K}_{\pi}}.$$
Since $\Phi_{\varphi}$ is a non-degenerate $\varphi$-map, the
above equality guarantees the existence of a unique isometry
$W_{\varphi}:\mathcal{H}_{\varphi}\to\mathcal{K}_{\pi}$ such that
$W_{\varphi}\Phi_{\varphi}(x)=\Psi_{\pi}(x)V$ satisfies for all
$x\in\mathcal{E}.$ Since $\Phi_{\varphi}$ and $\Psi_{\pi}$ are
non-degenerate continuous linear maps and $(\pi,\mathcal{K},V)$ is
a minimal Stinespring  dilation for $\varphi$,
\begin{equation*}\begin{split}
W_{\varphi}(\mathcal{H}_{\varphi})&=W_{\varphi}([\Phi_{\varphi}(\mathcal{E})\mathcal{H}])=[W_{\varphi}\Phi_{\varphi}(\mathcal{E})\mathcal{H}]
=[\Psi_{\pi}(\mathcal{E})V\mathcal{H}]\\&=
[\Psi_{\pi}(\mathcal{E})\pi(\mathcal{A})V\mathcal{H}]=[\Psi_{\pi}(\mathcal{E})[\pi(\mathcal{A})V\mathcal{H}]]
=[\Psi_{\pi}(\mathcal{E})\mathcal{K}]=\mathcal{K}_{\pi}
\end{split}\end{equation*}
so $W_{\varphi}$ is a unitary operator with the desired property.

 $\hspace{15cm} \square$

\vspace{0.5 cm}

\blue{\textbf{Proof of Theorem \ref{t2}:}}

\red{(i)} Similar to the proof of Lemma \ref{l1},
 for every
$x_{1},...,x_{n}\in\mathcal{E}$ and
$h_{1},...,h_{n}\in\mathcal{H}$ we have
$$\|\sum_{i=1}^{n}\Phi(x_{i})h_{i}\|_{\mathcal{H'}}=
\|\sum_{i=1}^{n}\Phi_{\varphi}(x_{i})h_{i}\|_{\mathcal{H}_{\varphi}}.$$

Thus there is an (onto) isometry
$S_{\Phi}:\mathcal{H}_\varphi\to\mathcal{H}'$ such that $S_{\Phi}
\Phi_\varphi(x)h=\Phi(x)h$ for every $x\in\mathcal{E}$ and
$h\in\mathcal{H}_{\Phi}.$
 Then
$\Phi(x)=S_{\Phi}W_\varphi^*\Psi_{\pi}(x)V$ for every
$x\in\mathcal{E}.$
 Put
$W:=W_\varphi S_{\Phi}^*,$ since $W_\varphi$ is a unitary and
$S_{\Phi}$ is an isometry, $W$ is a coisometry and
$\Phi(x)=W^*\Psi_{\pi}(x)V$ and $S_{\Phi} \Phi_\varphi(x)=\Phi(x)$
for every $x\in\mathcal{E}.$

For non-degenerate case, we have
$[\Phi(\mathcal{E})\mathcal{H}]=\mathcal{H'}$. Hence, isometry
$S_{\Phi}$ is onto and so it is unitary. Consequently, $W$ is a
unitary operator, too.

Conversely,  each of the equations
$\Phi(\cdot)=W^*\Psi_{\pi}(\cdot)V$ when $W$ is a coisometry  and
$\Phi(\cdot)=S_\Phi\Phi_\varphi(\cdot)$ when $S_\Phi$ is an
isometry, imply that $\Phi$ is a $\varphi$-map.

\red{(ii)} Let $\Phi$ be a (non-degenerate) completely
semi-$\varphi$-map. For every $x_1,...,x_n\in\mathcal{E}$, we have
  $$[\Phi(x_i)^*\Phi(x_j)]_{i,j}\leq [\varphi(\langle x_i,x_j\rangle)]_{i,j}.$$
 Consequently, for every $x_1,...,x_n\in\mathcal{E}$ and $h_1,...,h_n\in\mathcal{H}$ we have
 $$\|\sum_{i=1}^n\Phi(x_i)h_i\|^2\leq\sum_{i=1}^n\sum_{j=1}^n\langle\varphi(\langle x_j,x_i \rangle)h_i,h_j\rangle=
 \|\sum_{i=1}^n\Phi_{\varphi}(x_i)h_i\|^2.$$
Thus there is a (dense range) contractive linear operator
$S:\mathcal{H}_\varphi\to\mathcal{H}'$ such that
$S\Phi_\varphi(x)=\Phi(x)$ for every $x\in\mathcal{E}.$
 Therefore
$\Phi(x)=SW_\varphi^*\Psi_{\pi}(x)V$ for every $x\in\mathcal{E}.$
 Put
$W:=W_\varphi S^*,$ since $W_\varphi$ is a unitary and $S$ is a
(dense range) contractive operator, $W$ is a (injective)
contraction, and $\Phi(x)=W^*\Psi_{\pi}(x)V$ for every
$x\in\mathcal{E}.$

Conversely, when $W:\mathcal{H}_\varphi\to\mathcal{H}'$ is a
contraction, then the equation $\Phi(\cdot)=W^*\Psi(\cdot)V$
implies that $\Phi$ is a completely semi-$\varphi$-map.

\red{(iii)} It follows from (i).

\red{(iv)} For prove this, it is sufficient to set $W:=S_\Psi
W_\varphi S_{\Phi}^*.$

$\hspace{15cm} \square$

\blue{\textbf{Proof of Corollary \ref{c1}:}}

\red{$(i \Leftrightarrow ii)$}: By Theorem 2.1.

\red{($ii \Rightarrow iii$)}:  Let $\psi:=\|S\|^2 \varphi$.
Let $[x_{ij}] \in M_n(\mathcal{E})$. Thus
\begin{align*}
\Phi_{n}([x_{ij}])^*
\Phi_{n}([x_{ij}]) &= [\Gamma(x_{ji})^*S^*][S\Gamma(x_{ij})]=
 [\Gamma(x_{ji})^*]diag(S^*, \cdot\cdot\cdot, S^*) diag(S,
\cdot\cdot\cdot, S)[\Gamma(x_{ij})] \\&
 \leq ||S||^2 [\Gamma(x_{ji})^*][\Gamma(x_{ij})]
 = ||S||^2 \varphi_n(\langle
[x_{ij}],[x_{ij}]\rangle)
\end{align*}
and then $\Phi$ is completely semi-$\psi$-map.

 \red{($iii \Rightarrow iv$)}: By part $(ii)$ of Theorem \ref{t2}.

 \red{($iv \Rightarrow v$)}: See \cite[Theorem 4.2]{ABN}

 \red{($v \Rightarrow i$)}: As $\Phi$ is $1-2$ corner of some completely positive mapping on the linking $C^*$-algebra  $\mathcal{L}(\mathcal{E})$, then $\Phi$ is a completely bounded map.

 \subsection*{Acknowledgment}
The research of the first author was in part supported by a grant
from IPM (No. 94470046).


\begin{thebibliography}{99}
%
%


\bibitem{Arambašic}  L. Aramba$\check{s}$i$\acute{c}$, {\it Irreducible representations of Hilbert C*-modules}, Math. Proc. R. Ir. Acad. {2} (2005),  11-24.

\bibitem{Ars} N. Aronszajn, {\it Theory of Reproducing Kernels},  Tran. Amer. Math. Soc.  \textbf{68}  (3)  (1950), 337-404.

\bibitem{Asadi} M. B. Asadi, {\it Stinespring's theorem for Hilbert $C^*$-modules}, J. Operator Theory \textbf{62} (2) (2008), 235-238.

\bibitem{ABN} M. B. Asadi, R. Behmani, A. R. Medghalchi, H. Nikpey, {\it Completely semi-$\varphi$-maps}, arXiv:1608.00188.

\bibitem{BRS} B. V. R. Bhat, G. Ramesh and  K. Sumesh, {\it Stinespring's theorem for maps on Hilbert $C^*$-modules},
 J. Operator Theory \textbf{68} (2012), 173-178.

\bibitem{Dey-Tervidi} S. Dey and H. Trivedi, {\it $\mathfrak{K}$-families and CPD-H-extendable families}, arXive:1409.3655v1.

\bibitem{Gs} D. Goswami and K. B. Sinha, \textit{Quantum Stochastic Processes and  Noncommutative Geometry},
 Camb. Tra.  Math., Vol \textbf{169}, Camb. Uni. Press, 2007.

\bibitem{kol}A. N. Kolmogorov, {\it Stationary sequences in Hilbert space}, Bull. Math. Univ.
Moscow \textbf{2} (1941), 1–40.

\bibitem{Heo-Ji} J. Heo and Un C. Ji, {\it Quantum stochastic processes for maps on Hilbert $C^*$-modules}, J. Math. Phys. \textbf{52} (2011), 053501.

\bibitem{L} E. C. Lance, \textit{Hilbert $C^*$-modules}, Lond. Math. Soc. Lec. Note Ser., vol. \textbf{210}, Camb. Uni. Press, 1995.

\bibitem{J1} M. Joita, {\it Covariant version of the Stinespring type theorem for Hilbert $C^*$-modules},
Cent. Eur. J. Math. \textbf{9} (4)(2011), 803-813.

\bibitem{J2} M. Joita, {\it Comparison of completely positive maps on Hilbert $C^*$-modules}, J. Math. Anal. Appl.  \textbf{393} (2012), 644-650.

\bibitem{Pashke} W. L. Paschke, {\it Inner Product Modules Over $B^*$-Algebras},
  Tran. Amer. Math. Soc. \textbf{182} (1973), 443-468.

\bibitem{Pa1} V. Paulsen, \textit{\it An Introduction to the theory of Reproducing Kernel Hilbert Spaces}, www.math.uh.edu/vern, 2009.

\bibitem{Pau} V. Paulsen, \textit{Completely Bounded Maps and Operator Algebras},
Camb. Stud. adva. Math., vol. \textbf{78}, Camb. Uni. Press, 2002.


\bibitem{skeide} M. Skeide, {\it Factorization of maps between Hilbert $C^*$-modules}, J. Operator Theory \textbf{68} (2012), 543-547.

\bibitem{Skeide-Sumesh} M. Skeide, K. Sumesh, {\it CP-H-Extendable maps between Hilbert Modules and CPH-semigroups}, J. Math. Anal. Appl. \textbf{414} (2014), 886-913.


\bibitem{St} W. F. Stinespring, {\it Positive functions on $C^*$-algebras}, Proc. Amer. Math.
Soc. \textbf{6} (1955), 211-216.

\end{thebibliography}
\end{document}